\documentclass[12pt]{article}
\usepackage{graphicx}
\makeatletter
\let\normalequation=\equation
\def\equation{\@ifnextchar[{\subequation}{\normalequation}}
\def\subequation[#1]#2{\@ifundefined{r@#1}%
  {\def\theequation{\bf ??#2}\@warning
    {Reference `#1' on page \thepage \space
     undefined}}{\edef\@tempa{\@nameuse{r@#1}}%
    \edef\theequation{\expandafter\@car\@tempa \@nil#2}}%
  \let\@currentlabel\theequation $$}
\makeatother

\begin{document}

\begin{center}
{\large\bf Global Exponential Stability of Delayed Periodic Dynamical Systems}
\footnote{It is supported
  by National Science Foundation of China 60074005 and
  60374018}
\\[0.2in]
\begin{center}

\end{center}

\end{center}

\begin{center}
 Yanxu Zheng, Tianping Chen\footnote{These authors are with Lab. of
Nonlinear Mathematics Science, Institute of Mathematics, Fudan
University, Shanghai, 200433, P.R.China.\\
\indent ~~Corresponding author: Tianping Chen.
Email:tchen@fudan.edu.cn}
\end{center}

\section{Abstract}
In this paper, we discuss delayed periodic dynamical systems,
compare capability of criteria of global exponential stability in
terms of various $L^{p}$ ($1\le p<\infty$) norms. A general
approach to investigate global exponential stability in terms of
various $L^{p}$ ($1\le p<\infty$) norms is given.
Sufficient conditions ensuring global exponential stability
are given, too. Comparisons of various stability criteria are given.
More importantly, it is pointed out that sufficient conditions in
terms of $L^{1}$ norm are enough and easy to implement in practice.

\bigskip
{\bf Key words:} Neural Networks, Stability, Global exponential
convergence.

\pagestyle{plain} \pagenumbering{arabic}


\section{\bf Introduction}

Recurrently connected neural networks, sometimes called
Grossberg-Hopfield neural networks, have been extensively studied
in past years and found many applications in different areas.
However, many applications heavily depend on the dynamic behaviors
of the networks. Therefore, analysis of these dynamic behaviors is
a necessary step toward practical design of these neural networks.

A recurrently connected neural network is described by the
following differential equations:
\begin{equation}
\frac{du_{i}(t)}{dt}=-d_{i}u_{i}(t)+\sum_{j}a_{ij}g_{j}(u_{j}(t))+I_{i}
\quad i=1,\cdots,n \label{Hopfield}
\end{equation}
where $g_{j}(x)$ are activation functions, $d_{i}$, $a_{ij}$ are
constants and $I_{i}$ are constant inputs. In practice, however,
the interconnections are  generally asynchronous. Therefore, one
often needs to investigate the following delayed dynamical
systems:
\begin{equation}
\frac{du_{i}(t)}{dt}=-d_{i}u_{i}(t)+\sum_{j}a_{ij}g_{j}(u_{j}(t))
+\sum_{j}b_{ij}f_{j}(u_{j}(t-\tau_{ij}))+I_{i} \quad i=1,\cdots,n
\label{delay}
\end{equation}
 where activation functions $g_{j}$ and $f_{j}$ satisfy
certain defining conditions,  and
$a_{ij}$, $b_{ij}$, $I_{i}$ are constants.

All these neural networks and their corresponding delayed systems
have been extensively studied, and there are many papers
considering their global stability in the literature, for example,
see \cite{A,Cao,Chen1,Cohen,Forti,Hirsh,Liang,Liao,Roska,Roska1,Zhang,Zhang1}
and others.

However, the interconnection weights $a_{ij}$, $b_{ij}$,
self-inhibition $d_{i}$ and inputs $I_{i}$ should be variable with
time, often periodically. Therefore, we need to discuss following
dynamical systems and their periodic limits.
\begin{equation}
\frac{du_i}{dt}=-d_{i}(t)u_{i}(t)+\sum_{j=1}^{n}a_{ij}(t)g_j(u_j(t))
+\sum_{j=1}^{n}b_{ij}(t)f_{j}(u_{j}(t-\tau_{ij}))+I_i(t),\quad i=1,2,\ldots,n
\label{periodic}
\end{equation}
where    $d_{i}(t)\geq d_{i}>0$, $a_{ij}(t), b_{ij}(t), I_i(t):
\mathbf{R}^{+}\rightarrow \mathbf{R}$ are continuously periodic functions with period
$\omega>0$, i.e., $d_{i}(t+\omega)=d_{i}(t)$, $a_{ij}(t)=a_{ij}(t+\omega)$,
$b_{ij}(t)=b_{ij}(t+\omega)$, $I_{i}(t)=I_{i}(t+\omega)$ for all $t>0$ and
$i,j=1,2,\ldots,n$.

The initial conditions are
\begin{equation}
u_{i}(s)=\phi_{i}(s) \quad  for \quad s\in[-\tau,0], \quad where \quad \tau=\max_{1\le
i,j\le n}\tau_{ij}
\end{equation}
where $\phi_{i}\in C([-\tau,0]), i=1,\cdots,n$.

We also denote
$|a_{ij}^{*}|=\sup_{\{0<t\le \omega \}}|a_{ij}(t)|<\infty,$
$|b_{ij}^{*}|=\sup_{\{0<t\le \omega \}}|b_{ij}(t)|<\infty,$
$|I_{i}^{*}|=\sup_{\{0<t\le \omega \}}|I_{i}(t)|<\infty,$  where
$i,j=1,\cdots,n$.

\noindent There are also several papers discussing periodic
dynamical systems and their periodic solutions and its stability.
For example, see \cite{B,Cao1,Chen,Mo,Go,Xiang,Zhang,Zhou} and
others.

\noindent In this paper, our main concerns are comparisons of the
capability of various criteria in proving  the existence of
periodic solution and its stability. For this purpose, we first
discuss the existence of periodic solution and its stability in
terms of $L^{p}$ ($1\le p<\infty$) norms. We do not assume that
the activation functions are bounded. We also do not use existing
complicated theory (for example, topological degree theory, fixed
point theorem and so on), as needed in most other papers, to prove
the existence of periodic solution. Instead, we propose a general
and very concise approach.
By this approach, we prove exponential convergence directly.
The existence of the periodic solution is a direct consequence of
the exponential convergence.

\noindent This paper is organized in the following way: In section
3, some preliminaries, including several definitions and Young
Inequality as lemma, are given. In section 4, we prove the global
stability in terms of $L^{p}$ norms. Main comparison results are
given in section 5. In this section we point out that the criteria
in terms of $L^{1}$ norm are enough. We conclude the paper in
section 6.

\section{Some preliminaries}

\noindent{\bf Definition 1}\quad Class $H\{G_{1},\cdots,G_{n}\}$ of functions: Let
${G}=diag[{G}_{1},\cdots,{G}_{n}]$, where ${G}_{i}>0$, $i=1,\cdots,n$.
$g(x)=(g_{1}(x),\cdots,g_{n}(x))^{T}$ is said to belong to $H\{G_{1},\cdots,G_{n}\}$,
if the functions $g_{i}(x)$, $i=1,\cdots,n$ satisfy
$\frac{|g_{i}(x+u)-g_{i}(x)|}{|u|}\le G_{i}$.

\noindent{\bf Definition 2}\quad A vector $v=[v_{1},\cdots,v_{n}]^{T}>0$,
if and only if every $v_{i}>0$, $i=1,\cdots,n$.

\noindent{\bf Definition 3}\quad A real $n\times n$ matrix $C=(c_{ij})$ is said to be
an M-matrix if $c_{ij}\le 0$, $i,j=1,\cdots,n$, $j\neq i$ and all successive principal
minors of $C$ are positive.

\noindent{\bf Definition 4} Throughout this paper, we use the
following norm
\begin{equation}
||u(\cdot)||_{\{\xi,p\}}=\bigg[\frac{1}{n}
\sum_{i=1}^{n}\xi_{i}|u_{i}(\cdot)|^{p}\bigg]^{\frac{1}{p}}
\end{equation}
where $u(\cdot)=[u_{1}(\cdot),\cdots,u_{n}(\cdot)]^{T}$.

\noindent{\bf Lemma 1} (Young Inequality)\quad
If $a>0,\ b>0$, then
\begin{equation}
ab\leq\frac{(a\varepsilon)^{p}}{p}+\frac{(b\varepsilon^{-1})^{q}}{q}
\end{equation}
where
$\varepsilon >0;\ p,q>1;\ \frac{1}{p}+\frac{1}{q}=1$. The equality holds if and only if
\begin{equation}
(a\epsilon)^{p}=(b\epsilon^{-1})^{q}
\end{equation}

\noindent{\bf Proof}\quad
Let $f(x)=e^{x}$.
Obviously, $f''(x)>0$. Thus
$$f(\alpha x+\beta y)\leq \alpha f(x) + \beta f(y) \quad where \quad \alpha>0,\beta>0,\alpha +\beta=1$$
that's
\begin{equation}
e^{\alpha x}e^{\beta y}\leq \alpha e^{x}+ \beta e^{y}
\end{equation}
where equality  holds if and only if $x=y$.

Let $\alpha=\frac{1}{p},\ \ \beta=\frac{1}{q},\ \ x=\frac{1}{\alpha}\ln{a\varepsilon},\ \
y=\frac{1}{\beta}\ln{b\varepsilon^{-1}}$,
we obtain
\begin{equation}
ab\leq\frac{(a\varepsilon)^{p}}{p}+\frac{(b\varepsilon^{-1})^{q}}{q}
\end{equation}
and the equality holds if and only if
\begin{equation}
(a\epsilon)^{p}=(b\epsilon^{-1})^{q}
\end{equation}

\section{Main Results}
In this section, we discuss the existence of periodic solution and its stability.
We propose a general and concise approach and give several theorems on the existence of
periodic solution and its exponential stability.

{\bf Theorem 1}\quad Suppose that $1\le p<\infty$,
$g(x)=(g_{1}(x),\cdots,g_{n}(x))^{T}\in H\{G_{1},\cdots,G_{n}\}$
and $f(x)= (f_{1}(x),\cdots,f_{n}(x))^{T}\in
H\{F_{1},\cdots,F_{n}\}$.
 If there are real constants  $\epsilon>0$, $\xi_{i}>0$,
$\alpha_{ij}$, $\beta_{ij}$,  $i,j=1,2\cdots,n$,
 such that
\begin{eqnarray}
&&\nonumber (-d_{i}+\epsilon)\xi_{i}+G_{i}\bigg[\xi_{i}|a_{ii}^{*}|
+\frac{1}{p}\sum\limits_{j\ne i}\xi_{j}|a_{ji}^{*}|^{\alpha_{ji}p}\bigg]
+\frac{1}{q}\xi_{i}\sum\limits_{j\ne
i}G_{j}|a_{ij}^{*}|^{(1-\alpha_{ij})q}\\
&+&\frac{1}{p}F_{i}\sum\limits_{j=1}^{n}\xi_{j}|b^{*}_{ji}|^{\beta_{ji}p}e^{\epsilon
\tau_{ji}}+
\frac{1}{q}\xi_{i}\sum\limits_{j=1}^{n}F_{j}|b^{*}_{ij}|^{(1-\beta_{ij})q}
e^{\epsilon \tau_{ij}} \le 0 \label{Th1-modified}
\end{eqnarray}
In particular,  $\alpha_{ij}=\beta_{ij}=\frac{1}{p}$,
\begin{eqnarray}
&&\nonumber (-d_{i}+\epsilon)\xi_{i}+G_{i}\bigg[\xi_{i}|a_{ii}^{*}|
+\frac{1}{p}\sum\limits_{j\ne i}\xi_{j}|a_{ji}^{*}|\bigg]
+\frac{1}{q}\xi_{i}\sum\limits_{j\ne
i}G_{j}|a_{ij}^{*}|\\
&+&\frac{1}{p}F_{i}\sum\limits_{j=1}^{n}\xi_{j}|b^{*}_{ji}|e^{\epsilon
\tau_{ji}}+
\frac{1}{q}\xi_{i}\sum\limits_{j=1}^{n}F_{j}|b^{*}_{ij}|
e^{\epsilon \tau_{ij}} \le 0 \label{Th1}
\end{eqnarray}
Then  the dynamical system (\ref{periodic}) has a unique periodic
solution $v(t)=[v_{1}(t),\cdots,v_{n}(t)]$ and, for any solution
$u(t)=[u_{1}(t),\cdots,u_{n}(t)]$ of (\ref{periodic}),
\begin{equation}
|u_{i}(t+j\omega)-v_{i}(t)|=O(e^{-j\epsilon\omega }), \qquad
i=1,\cdots,n
\end{equation}

\noindent{\bf Proof}\quad \quad
 Let
$\bar{u}_{i}(t)=u_{i}(t+\omega)-u_{i}(t),$
$\bar{g}_{i}(u_{i}(t))=g_{i}(u_{i}(t+\omega))-g_{i}(u_{i}(t))$,
$\bar{f}_{i}(u_{i}(t))=f_{i}(u_{i}(t+\omega))-f_{i}(u_{i}(t))$,
and $w_{i}(t)=e^{\epsilon t}\bar{u}_{i}(t)$, $ i=1,2,\cdots,n.$

Defining a Lyapunov function by
\begin{equation}
L_{1}(t)=\sum_{i=1}^{n}\xi_{i}|w_{i}(t)|^{p}
+p\sum_{i,j=1}^{n}\xi_{i}F_{j}|b_{ij}^{*}|e^{\epsilon \tau_{ij}}
\int_{t-\tau_{ij}}^{t}|w_{j}(y)|^pdy
\end{equation}
and differentiating it, we have
\begin{eqnarray}
&&\dot{L}_{1}(t)=p\sum\limits_{i=1}^{n}\xi_{i}|w_{i}(t)|^{p-1}e^{\epsilon t
}sign(w_{i}(t))
\bigg[-d_{i}(t)\bar{u}_{i}(t)+\epsilon \bar{u}_{i}(t)\nonumber\\
&+&\sum\limits_{j=1}^{n}a_{ij}(t)\bar{g}_{j}(u_{j}(t))
+\sum\limits_{j=1}^{n}b_{ij}(t)\bar{f}_{j}(u_{j}(t-\tau_{ij}))\bigg]\nonumber\\
&+&p\sum_{i,j=1}^{n}\xi_{i}F_{j}|b_{ij}^{*}|e^{\epsilon
\tau_{ij}}\bigg[|w_{j}(t)|^{p}
-|w_{j}(t-\tau_{ij})|^p\bigg]\\
&\le&
p\sum\limits_{i=1}^{n}\xi_{i}\bigg[-(d_{i}-\epsilon)|w_{i}(t)|^{p}
+a_{ii}(t)e^{\epsilon t}|\bar{g}_{i}(u_{i}(t))||w_{i}(t)|^{p-1}\nonumber\\
&+&\sum\limits_{j\ne i}e^{\epsilon
t}|a_{ij}^{*}||w_{i}(t)|^{p-1}|\bar{g}_{j}(u_{j}(t))|
+\sum\limits_{j=1}^{n} |b_{ij}^{*}|e^{\epsilon
t}|\bar{f}_{j}(u_{j}(t-\tau_{ij}))||w_{i}(t)|^{p-1}
\bigg]\nonumber\\
&+&p\sum_{i,j=1}^{n}\xi_{i}F_{j}|b_{ij}^{*}|e^{\epsilon \tau_{ij}}
\bigg[|w_{j}(t)|^{p} -|w_{j}(t-\tau_{ij})|^p\bigg]
\end{eqnarray}
By Young inequality and
\begin{equation}
e^{\epsilon t}|\bar{f}_{j}(u_{j}(t-\tau_{ij}))|\le
F_{j}|w_{j}(t-\tau_{ij})|e^{\epsilon \tau_{ij}}
\end{equation}
\begin{equation}
e^{\epsilon t}|\bar{g}_{j}(u_{j}(t))|\le G_{j}|w_{j}(t)|
\end{equation}
 we have
\begin{eqnarray}
\dot{L}_{1}(t)&\le& p
\sum\limits_{i=1}^{n}\bigg\{-(d_{i}-\epsilon)\xi_{i}+G_{i}\bigg[\xi_{i}|a_{ii}^{*}|
+\frac{1}{p}\sum\limits_{j\ne i}|a_{ji}^{*}|^{\alpha_{ji}p}\xi_{j}\bigg]
+\frac{1}{q}\xi_{i}\sum\limits_{j\ne
i}G_{j}|a_{ij}^{*}|^{(1-\alpha_{ij})q}\nonumber\\
&+&\sum\limits_{j=1}^{n}\bigg[\frac{1}{p}\xi_{j}F_{i}|b_{ji}^{*}|^{\beta_{ji}p}
e^{\epsilon \tau_{ji}}
+\frac{1}{q}\xi_{i}F_{j}|b_{ij}^{*}|^{(1-\beta_{ij})q}e^{\epsilon \tau_{ij}}
\bigg]\bigg\} |w_{i}(t)|^{p}
\nonumber\\
&\le& 0
\end{eqnarray}
\ Therefore, $L_{1}(t)$ is bounded, which implies
\begin{equation}
\sum_{i=1}^{n}\xi_{i}|e^{\epsilon t}\bar{u}_{i}(t)|^{p}
\end{equation}
is bounded and
\begin{eqnarray}
|u_{i}(t+\omega)-u_{i}(t)| =O(e^{-\epsilon t}), \qquad
i=1,\cdots,n
 \label{LL1}
\end{eqnarray}
 Now, define a function $v(t)=[v_{1}(t),\cdots,v_{n}(t)]^{T}$ by
 \begin{eqnarray*}
v_{i}(t)=\lim_{j\rightarrow \infty}u_{i}(t+j\omega)
\end{eqnarray*}
Because of
\begin{eqnarray*}
u_{i}(t+j\omega)=u_{i}(t)+\sum_{k=1}^{j}\bigg\{u_{i}(t+k\omega)-u_{i}(t+(k-1)\omega)\bigg\}
\end{eqnarray*}
and (\ref{LL1}), $v(t)$ is well defined and is a periodic function with period
$\omega$. Moreover,

If $u(t)$, $v(t)$ are two solutions. By similar method used
before, it is easy to prove
 \begin{eqnarray}
|u_{i}(t+j\omega)-v_{i}(t+j\omega)| =O(e^{-j\epsilon \omega})
\quad when \quad j\rightarrow \infty
\end{eqnarray}
which means the limit solution is unique. Theorem 1 is proved
completely.

{\bf Corollary 1}\quad Suppose that $1\le p<\infty$,
$g(x)=(g_{1}(x),\cdots,g_{n}(x))^{T}\in H\{G_{1},\cdots,G_{n}\}$
and $f(x)= (f_{1}(x),\cdots,f_{n}(x))^{T}\in
H\{F_{1},\cdots,F_{n}\}$.
 If there are real constants  $\xi_{i}>0$,
$\alpha_{ij}$, $\beta_{ij}$,   $i,j=1,2\cdots,n$,
 such that
\begin{eqnarray}
&&\nonumber -d_{i}\xi_{i}+G_{i}\bigg[\xi_{i}|a_{ii}^{*}| +\frac{1}{p}\sum\limits_{j\ne
i}\xi_{j}|a_{ji}^{*}|^{\alpha_{ji}p}\bigg] +\frac{1}{q}\xi_{i}\sum\limits_{j\ne
i}G_{j}|a_{ij}^{*}|^{(1-\alpha_{ij})q}\\
&+&\frac{1}{p}F_{i}\sum\limits_{j=1}^{n}\xi_{j}|b^{*}_{ji}|^{\beta_{ji}p}+
\frac{1}{q}\xi_{i}\sum\limits_{j=1}^{n}F_{j}|b^{*}_{ij}|^{(1-\beta_{ij})q}
 < 0 \label{Co1}
\end{eqnarray}
In particular,
\begin{eqnarray}
&&\nonumber -d_{i}\xi_{i}+G_{i}\bigg[\xi_{i}|a_{ii}^{*}| +\frac{1}{p}\sum\limits_{j\ne
i}\xi_{j}|a_{ji}^{*}|\bigg] +\frac{1}{q}\xi_{i}\sum\limits_{j\ne
i}G_{j}|a_{ij}^{*}|\\
&+&\frac{1}{p}F_{i}\sum\limits_{j=1}^{n}\xi_{j}|b^{*}_{ji}|+
\frac{1}{q}\xi_{i}\sum\limits_{j=1}^{n}F_{j}|b^{*}_{ij}|
 < 0
\label{cc}
\end{eqnarray}
Then  the dynamical system (\ref{periodic}) has a unique periodic
solution $v(t)=[v_{1}(t),\cdots,v_{n}(t)]^{T}$ and, there is
$\epsilon>0$ such that for any solution
$u(t)=[u_{1}(t),\cdots,u_{n}(t)]^{T}$ of (\ref{periodic}), we have
 \begin{eqnarray}
|u_{i}(t+j\omega)-v_{i}(t)| =O(e^{-j\epsilon \omega}) \quad when
\quad j\rightarrow \infty
\end{eqnarray}

In fact, under the assumptions given in Corollary 1, we can find $\epsilon$ such that
(\ref{Th1-modified}) or (\ref{Th1}) is satisfied.

The case $p=1$ is the most interesting. In Theorem 1, let $p=1$ and
$\alpha_{ij}=\beta_{ij}=1$, we have following

{\bf Theorem 2}\quad
Suppose that $1\le p<\infty$,
$g(x)=(g_{1}(x),\cdots,g_{n}(x))^{T}\in H\{G_{1},\cdots,G_{n}\}$
and $f(x)= (f_{1}(x),\cdots,f_{n}(x))^{T}\in
H\{F_{1},\cdots,F_{n}\}$.
If there are positive constants $\epsilon$,
$\theta_{i}$,  $i=1,2\cdots,n$,
 such that
\begin{eqnarray}
&& (-d_{i}+\epsilon)\theta_{i}+G_{i}\bigg[\theta_{i}|a_{ii}^{*}|
+\sum\limits_{j\ne i}\theta_{j}|a_{ji}^{*}|\bigg]
+F_{i}\sum\limits_{j=1}^{n}\theta_{j}|b^{*}_{ji}|e^{\epsilon
\tau_{ji}} \le 0 \label{Th2}
\end{eqnarray}
Then  the dynamical system (\ref{periodic}) has a unique periodic
solution $v(t)=[v_{1}(t),\cdots,v_{n}(t)]^{T}$ and, for any
solution $u(t)=[u_{1}(t),\cdots,u_{n}(t)]^{T}$ of
(\ref{periodic}),
\begin{equation}
|u_{i}(t+j\omega)-v_{i}(t)|=O(e^{-j\epsilon\omega }), \qquad
i=1,\cdots,n
\end{equation}
If
\begin{eqnarray}
&& -d_{i}\theta_{i}+G_{i}\bigg[\theta_{i}|a_{ii}^{*}|
+\sum\limits_{j\ne i}\theta_{j}|a_{ji}^{*}|\bigg]
+F_{i}\sum\limits_{j=1}^{n}\theta_{j}|b^{*}_{ji}| < 0 \label{Th21}
\end{eqnarray}
Then  the dynamical system (\ref{periodic}) has a unique periodic
solution $v(t)=[v_{1}(t),\cdots,v_{n}(t)]^{T}$ and, there is
$\epsilon>0$ such that for any solution
$u(t)=[u_{1}(t),\cdots,u_{n}(t)]^{T}$ of (\ref{periodic}), we have
 \begin{eqnarray}
|u_{i}(t+j\omega)-v_{i}(t)| =O(e^{-j\epsilon \omega}), \qquad
i=1,\cdots,n
\end{eqnarray}

Theorem 1 and Theorem 2 apply to the model (\ref{delay}), too. If
we consider constants $d_{i}$, $a_{ij}$, $b_{ij}$ $I_{i}$ as
periodic function with any period. Then, the limit $v(t)$ is also
a periodic function with any period and thus is a constant vector
$v^{*}=[v_{1}^{*},\cdots,v_{n}^{*}]^{*}$. Therefore, we have

{\bf Theorem 3}\quad Suppose that $1\le p<\infty$,
$g(x)=(g_{1}(x),\cdots,g_{n}(x))^{T}\in H\{G_{1},\cdots,G_{n}\}$
and $f(x)= (f_{1}(x),\cdots,f_{n}(x))^{T}\in
H\{F_{1},\cdots,F_{n}\}$.
 If there are real constants $\alpha_{ij}$,
$\beta_{ij}$, positive constants $\epsilon$, $\xi_{i}$, $\theta_{i}$,
  $i,j=1,2\cdots,n$,
 such that either one set on inequalities holds
\begin{eqnarray}
&&\nonumber
(-d_{i}+\epsilon)\xi_{i}+G_{i}\bigg[\xi_{i}|a_{ii}|
+\frac{1}{p}\sum\limits_{j\ne
i}\xi_{j}|a_{ji}|^{\alpha_{ji}p}\bigg]
+\frac{1}{q}\xi_{i}\sum\limits_{j\ne
i}G_{j}|a_{ij}|^{(1-\alpha_{ij})q}\\
&+&\frac{1}{p}F_{i}\sum\limits_{j=1}^{n}\xi_{j}|b_{ji}|^{\beta_{ji}p}e^{\epsilon
\tau_{ji}}+ \frac{1}{q}\xi_{i}\sum\limits_{j=1}^{n}F_{j}|b_{ij}|^{(1-\beta_{ij})q}
e^{\epsilon \tau_{ij}} \le 0 \label{Th31}
\end{eqnarray}
\begin{eqnarray}
&&\nonumber
(-d_{i}+\epsilon)\xi_{i}+G_{i}\bigg[\xi_{i}|a_{ii}|
+\frac{1}{p}\sum\limits_{j\ne i}\xi_{j}|a_{ji}|\bigg]
+\frac{1}{q}\xi_{i}\sum\limits_{j\ne
i}G_{j}|a_{ij}|\\
&+&\frac{1}{p}F_{i}\sum\limits_{j=1}^{n}\xi_{j}|b_{ji}|e^{\epsilon \tau_{ji}}+
\frac{1}{q}\xi_{i}\sum\limits_{j=1}^{n}F_{j}|b_{ij}| e^{\epsilon \tau_{ij}} \le 0
\label{Th32}
\end{eqnarray}
\begin{eqnarray}
&& (-d_{i}+\epsilon)\theta_{i}+G_{i}\bigg[\theta_{i}|a_{ii}|
+\sum\limits_{j\ne i}\theta_{j}|a_{ji}|\bigg]
+F_{i}\sum\limits_{j=1}^{n}\theta_{j}|b_{ji}|e^{\epsilon
\tau_{ji}} \le 0 \label{Th33}
\end{eqnarray}
Then  the dynamical system (\ref{delay}) has a unique periodic
equilibrium $v^{*}=[v_{1}^{*},\cdots,v_{n}^{*}]^{T}$ and, for any
solution $u(t)=[u_{1}(t),\cdots,u_{n}(t)]^{T}$ of (\ref{delay}),
we have
\begin{equation}
|u_{i}(t)-v_{i}^{*}|=O(e^{-\epsilon t }), \qquad
i=1,\cdots,n
\end{equation}

{\bf Corollary 2}\quad Suppose that $1\le p<\infty$,
$g(x)=(g_{1}(x),\cdots,g_{n}(x))^{T}\in H\{G_{1},\cdots,G_{n}\}$
and $f(x)= (f_{1}(x),\cdots,f_{n}(x))^{T}\in
H\{F_{1},\cdots,F_{n}\}$.
 If there are real constants $\alpha_{ij}$,
$\beta_{ij}$, positive constants  $\xi_{i}$, $\theta_{i}$,
  $i,j=1,2\cdots,n$,
 such that either one set on inequalities holds
\begin{eqnarray}
&&\nonumber
-d_{i}\xi_{i}+G_{i}\bigg[\xi_{i}|a_{ii}|
+\frac{1}{p}\sum\limits_{j\ne
i}\xi_{j}|a_{ji}|^{\alpha_{ji}p}\bigg]
+\frac{1}{q}\xi_{i}\sum\limits_{j\ne
i}G_{j}|a_{ij}|^{(1-\alpha_{ij})q}\\
&+&\frac{1}{p}F_{i}\sum\limits_{j=1}^{n}\xi_{j}|b_{ji}|^{\beta_{ji}p}+
\frac{1}{q}\xi_{i}\sum\limits_{j=1}^{n}F_{j}|b_{ij}|^{(1-\beta_{ij})q}
 < 0 \label{Co21}
\end{eqnarray}
\begin{eqnarray}
&&\nonumber
-d_{i}\xi_{i}+G_{i}\bigg[\xi_{i}|a_{ii}|
+\frac{1}{p}\sum\limits_{j\ne i}\xi_{j}|a_{ji}|\bigg]
+\frac{1}{q}\xi_{i}\sum\limits_{j\ne
i}G_{j}|a_{ij}|\\
&+&\frac{1}{p}F_{i}\sum\limits_{j=1}^{n}\xi_{j}|b_{ji}|+
\frac{1}{q}\xi_{i}\sum\limits_{j=1}^{n}F_{j}|b_{ij}|
 < 0 \label{Co22}
\end{eqnarray}
\begin{eqnarray}
-d_{i}\theta_{i}+G_{i}\bigg[\theta_{i}|a_{ii}| +\sum\limits_{j\ne
i}\theta_{j}|a_{ji}|\bigg] +F_{i}\sum\limits_{j=1}^{n}\theta_{j}|b_{ji}| < 0
\label{Co23}
\end{eqnarray}
Then  the dynamical system (\ref{delay}) has a unique periodic
equilibrium $v^{*}=[v_{1}^{*},\cdots,v_{n}^{*}]^{T}$ and, there is
$\epsilon>0$ such that for any solution
$u(t)=[u_{1}(t),\cdots,u_{n}(t)]^{T}$ of (\ref{delay}), we have
 \begin{eqnarray}
|u_{i}(t)-v_{i}(t)| =O(e^{-\epsilon t}),  \qquad
i=1,\cdots,n
\end{eqnarray}

{\bf Remark 1}\quad In \cite{Cao,D}, under additional assumptions
that all $a_{ij}(t)$, $b_{ij}(t)$ and $d_{i}(t)$ are constants and
$g_{j}(x)$, $f_{j}(x)$ are bounded functions, it was proved in
\cite{Cao} that under (\ref{cc}) or in \cite{D} that under
(\ref{Th1-modified}), the dynamical system
\begin{equation}
\frac{du_i}{dt}=-d_{i}u_{i}(t)+\sum_{j=1}^{n}a_{ij}g_j(u_j(t))
+\sum_{j=1}^{n}b_{ij}f_{j}(u_{j}(t-\tau_{ij}))+I_i(t),\quad
i=1,2,\ldots,n \label{cao}
\end{equation}
has a unique periodic solution
$v(t)=[v_{1}(t),\cdots,v_{n}(t)]^{T}$ and,  any solution
$u(t)=[u_{1}(t),\cdots,u_{n}(t)]^{T}$ of (\ref{cao}) converges to
$v(t)$. It is easy to see that this is a special case of previous
Corollary 1. Moreover, in \cite{Cao}, the author did not address
the convergence rate. It should also be emphasized that parameters
$\alpha_{ij}$, $\beta_{ij}$ play key role in the Theorem 1 and
Corollary 1. The special one (\ref{cc}) is not the best criterion.
It is important to look for the best parameters $\alpha_{ij}$,
$\beta_{ij}$. We will discuss in next section.


\section{Comparisons}
\noindent In previous section and many existing papers, various
stability criteria are given. It is an important issue to answer
the question whether or not they are equivalent. If not, which one
is better. That is  we should compare capability of various
stability criteria.

\noindent In this section, we will address this issue by proving
several theorems, from which we assert that the conditions in terms
of $L^{1}$ norm [(\ref{Th2}) (\ref{Th21}) (\ref{Th33})
(\ref{Co23})] are enough in practice. The criteria with $L^{p}$
norm [(\ref{Th1-modified}) (\ref{Co1}) (\ref{Th31}) (\ref{Co21})]
are of no great significance and those, letting $\alpha_{ij}=1/p,\
\beta_{ij}=1/p\ $ [(\ref{Th1}) (\ref{cc}) (\ref{Th32})
(\ref{Co22})], are  less capable than conditions in terms of
$L^{1}$.

\noindent We prove the following result first:

\noindent{\bf Theorem 4} \quad Suppose $c_{i}>0,$  $G_{i}>0,$
$F_{i}>0$, $c_{ij}\geq 0$, $d_{ij}\geq 0$, $e_{ij}\geq 0$,
$i,j=1,\cdots,n$. If there exist $\xi_{i}>0,\ \alpha_{ij},\
\beta_{ij},\ p>1,\ q>1,\ \frac{1}{p}+\frac{1}{q}=1\ $ such that
the following inequalities
\begin{eqnarray}
&-&c_{i}\xi_{i}+\frac{1}{p}\sum\limits_{j\neq i}c_{ji}^{\alpha_{ji}p}G_{i}\xi_{j}
+\frac{1}{q}\sum\limits_{j\neq i}c_{ij}^{(1-\alpha_{ij})q}G_{j}\xi_{i}
+\frac{1}{p}\sum\limits_{j\neq i}d_{ji}^{\beta_{ji}p}e_{ji}F_{i}\xi_{j}\nonumber\\
&+&\frac{1}{q}\sum\limits_{j\neq
i}d_{ij}^{(1-\beta_{ij})q}e_{ij}F_{j}\xi_{i}\le 0, \qquad
i=1,\cdots,n\label{th41}
\end{eqnarray}
hold. Then, we can find constants $\theta_{i}>0\ (i=1,\cdots\,n)\
$ such that
\begin{eqnarray}
-c_{i}\theta_{i}+G_{i}\sum_{j\ne
i}\theta_{j}c_{ji}+F_{i}\sum_{j\ne i}\theta_{j}d_{ji}e_{ji}\le 0,
\qquad i=1,\cdots,n\label{th42}
\end{eqnarray}
hold.

\noindent{\bf Proof} \quad
If there exist $\xi_{i}>0,\ \alpha_{ij},\ \beta_{ij},\ p>1,\ q>1,\
\frac{1}{p}+\frac{1}{q}=1\ $ such that
(\ref{th41}) hold. Denote
$$
M=\{m_{ij}\}:\left\{
\begin{array}{ll}
m_{ii} =  -c_{i}+\frac{1}{q}\sum\limits_{j\neq
i}c_{ij}^{(1-\alpha_{ij})q}G_{j} +\frac{1}{q}\sum\limits_{j\neq
i}d_{ij}^{(1-\beta_{ij})q}e_{ij}F_{j}
& i=1,\cdots,n \\
m_{ij} =
\frac{1}{p}c_{ji}^{\alpha_{ji}p}G_{i}+\frac{1}{p}d_{ji}^{\beta_{ji}p}e_{ji}F_{i}
& i\neq j
\end{array}
\right.
$$

Then, (\ref{th41}) can be rewritten as $M\xi\le 0$. By the property of
M-matrices, there exist $\eta=(\eta_{1},\cdots,\eta_{n})^{T},\eta_{i}>0\ i=1,\cdots,n\ $
such that $M^{T}\eta\le 0$. That is
\begin{eqnarray}
&-&\bigg(c_{i}-\frac{1}{q}\sum\limits_{j\neq i}c_{ij}^{(1-\alpha_{ij})q}G_{j}
-\frac{1}{q}\sum\limits_{j\neq i}d_{ij}^{(1-\beta_{ij})q}e_{ij}F_{j}\bigg)\eta_{i}\nonumber\\
&+&\frac{1}{p}\sum\limits_{j\neq i}(c_{ij}^{\alpha_{ij}p}G_{j}
+d_{ij}^{\beta_{ij}p}e_{ij}F_{j})\eta_{j}\le 0\nonumber\\
\end{eqnarray}
which can be rewritten as
\begin{eqnarray}
&-&c_{i}\eta_{i}
+\sum\limits_{j\neq i}\bigg (\frac{1}{q}c_{ij}^{(1-\alpha_{ij})q}G_{j}\eta_{i}
+\frac{1}{p}c_{ij}^{\alpha_{ij}p}G_{j}\eta_{j}\bigg )\nonumber\\
&+&\sum\limits_{j\neq i}\bigg
(\frac{1}{q}d_{ij}^{(1-\beta_{ij})q}e_{ij}F_{j}\eta_{i}
                                +\frac{1}{p}d_{ij}^{\beta_{ij}p}e_{ij}F_{j}\eta_{j}\bigg )
\le 0\nonumber\\
\end{eqnarray}

By Lemma 1, we have
\begin{eqnarray}
\frac{1}{q}c_{ij}^{(1-\alpha_{ij})q}G_{j}\eta_{i}
    +\frac{1}{p}c_{ij}^{\alpha_{ij}p}G_{j}\eta_{j}
\geq c_{ij}G_{j}\eta_{i}^{\frac{1}{q}}\eta_{j}^{\frac{1}{p}}
\label{L1}
\end{eqnarray}
\begin{eqnarray}
\frac{1}{q}d_{ij}^{(1-\beta_{ij})q}e_{ij}F_{j}\eta_{i}
    +\frac{1}{p}d_{ij}^{\beta_{ij}p}e_{ij}F_{j}\eta_{j}
\geq d_{ij}e_{ij}F_{j}\eta_{i}^{\frac{1}{q}}\eta_{j}^{\frac{1}{p}}
\label{L2}
\end{eqnarray}
Therefore,
\begin{eqnarray}
-c_{i}\eta_{i} +\sum\limits_{j\neq
i}c_{ij}G_{j}\eta_{i}^{\frac{1}{q}}\eta_{j}^{\frac{1}{p}}
+\sum\limits_{j\neq
i}d_{ij}e_{ij}F_{j}\eta_{i}^{\frac{1}{q}}\eta_{j}^{\frac{1}{p}}
\le 0 \ \ \ i=1,\cdots,n\nonumber
\end{eqnarray}
Thus, let $\zeta_{i}=\eta_{i}^{\frac{1}{p}},\ (i=1,\cdots,n)$,
 we obtain
\begin{eqnarray*}
-c_{i}\zeta_{i}+\sum_{j\ne i}\zeta_{j}c_{ij}G_{j}+\sum_{j\ne
i}\zeta_{j}d_{ij}e_{ij}F_{j}\le 0, \qquad i=1,\cdots,n
\end{eqnarray*}
By the property of M-matrices, there exist $\theta_{i}$ s.t.
\begin{eqnarray*}
-c_{i}\theta_{i}+G_{i}\sum_{j\ne
i}\theta_{j}c_{ji}+F_{i}\sum_{j\ne i}\theta_{j}d_{ji}e_{ji}\le 0,
\qquad i=1,\cdots,n
\end{eqnarray*}
This completes the proof.

\noindent{\bf Theorem 5} \quad Suppose $ p>1,\ q>1,\
\frac{1}{p}+\frac{1}{q}=1\ $. $c_{i}>0,\ G_{i}>0,\ F_{i}>0\ \
(i,j=1,\cdots,n)$. If there exist $\theta_{i}>0\ (i=1,\cdots\,n)\
$ such that
\begin{eqnarray}
-c_{i}\theta_{i}+G_{i}\sum_{j\ne
i}\theta_{j}c_{ji}+F_{i}\sum_{j\ne i}\theta_{j}d_{ji}e_{ji}\le 0,
\qquad i=1,\cdots,n \label{th51}
\end{eqnarray}
hold. Then, we can find constants $\xi_{i}>0,\ \alpha^{*}_{ij},\
\beta^{*}_{ij},$ such that following inequalities
\begin{eqnarray}
&-&c_{i}\xi_{i}+\frac{1}{p}\sum\limits_{j\neq
i}c_{ji}^{\alpha^{*}_{ji}p}G_{i}\xi_{j}
+\frac{1}{q}\sum\limits_{j\neq
i}c_{ij}^{(1-\alpha^{*}_{ij})q}G_{j}\xi_{i}
+\frac{1}{p}\sum\limits_{j\neq i}d_{ji}^{\beta^{*}_{ji}p}e_{ji}F_{i}\xi_{j}\nonumber\\
&+&\frac{1}{q}\sum\limits_{j\neq
i}d_{ij}^{(1-\beta^{*}_{ij})q}e_{ij}F_{j}\xi_{i}\le 0, \qquad
i=1,\cdots,n\label{th52}
\end{eqnarray}
hold.

{\bf Proof}\quad If there exist $\theta_{i}>0\ (i=1,\cdots\,n)\ $
such that (\ref{th51}) hold then, by the property of M-matrices,
there exist $\zeta_{i}$ s.t.
\begin{eqnarray*}
-c_{i}\zeta_{i}+\sum_{j\ne i}\zeta_{j}c_{ij}G_{j}+\sum_{j\ne
i}\zeta_{j}d_{ij}e_{ij}F_{j}\le 0, \qquad i=1,\cdots,n
\end{eqnarray*}
Let $\eta_{i}=\zeta_{i}^{p}$, $i=1,\cdots,n,$ and
\begin{eqnarray}
\alpha^{*}_{ij}=\frac{1}{p}\bigg[1+\ln_{|c_{ij}|}\bigg(\eta_{i}^{\frac{1}{q}}
\eta_{j}^{-\frac{1}{q}}\bigg)\bigg] \quad
\beta^{*}_{ij}=\frac{1}{p}\bigg[1+\ln_{|d_{ij}|}\bigg(\eta_{i}^{\frac{1}{q}}
\eta_{j}^{-\frac{1}{q}}\bigg)\bigg]
\end{eqnarray}
(\ref{L1}), (\ref{L2}) turn to be equalities.
Therefore, we have
\begin{eqnarray*}
& &    -c_{i}\eta_{i}
        +\sum\limits_{j\neq i}\bigg (\frac{1}{q}c_{ij}^{(1-\alpha^{*}_{ij})q}G_{j}\eta_{i}
         +\frac{1}{p}c_{ij}^{\alpha^{*}_{ij}p}G_{j}\eta_{j}\bigg )\\
        &+&\sum\limits_{j\neq i}\bigg (\frac{1}{q}d_{ij}^{(1-\beta^{*}_{ij})q}e_{ij}F_{j}\eta_{i}
                                +\frac{1}{p}d_{ij}^{\beta^{*}_{ij}p}e_{ij}F_{j}\eta_{j}\bigg )\\
& = &   -c_{i}\eta_{i}
        +\sum\limits_{j\neq i}c_{ij}G_{j}\eta_{i}^{\frac{1}{q}}\eta_{j}^{\frac{1}{p}}
        +\sum\limits_{j\neq i}d_{ij}e_{ij}F_{j}\eta_{i}^{\frac{1}{q}}\eta_{j}^{\frac{1}{p}}\\
& = &   \eta_{i}^{\frac{1}{q}}\bigg ( -c_{i}\theta_{i}
                                    +\sum\limits_{j\neq i}c_{ij}G_{j}\theta_{j}
                                    +\sum\limits_{j\neq i}d_{ij}e_{ij}F_{j}\theta_{j} \bigg)
        \le 0\ \ \ \ i=1,\cdots,n
\end{eqnarray*}
which means $M^{T}\eta\le 0$. By M-matrices theory, there exist
$\xi=[\xi_{1},\cdots,\xi_{n}]^{T}>0$, such that $M\xi\le 0$. Thus,
(\ref{th52}) hold. This completes the proof.

Under that
$\frac{1}{p}|b_{ii}|^{\beta_{ii}p}+\frac{1}{q}|b_{ii}|^{(1-\beta_{ii}q)}\ge
|b_{ii}|$, letting
$c_{i}=d_{i}-\epsilon-G_{i}|a_{ii}|-F_{i}|b_{ii}|$,
$c_{ij}=|a_{ij}|$, $d_{ij}=|b_{ij}|$, $e_{ij}=e^{\epsilon
\tau_{ij}}$, we have

{\bf Theorem 6}\quad Suppose that $1\le p<\infty$.
 If there are positive constants $\epsilon$, $\xi_{i}$,
 real constants $\alpha_{ij}$, $\beta_{ij}$,
  $i,j=1,2\cdots,n$, such that either one of the following two sets of inequalities holds
\begin{eqnarray}
&&\nonumber
(-d_{i}+\epsilon)\xi_{i}+G_{i}\bigg[\xi_{i}|a_{ii}|
+\frac{1}{p}\sum\limits_{j\ne
i}\xi_{j}|a_{ji}|^{\alpha_{ji}p}\bigg]
+\frac{1}{q}\xi_{i}\sum\limits_{j\ne
i}G_{j}|a_{ij}|^{(1-\alpha_{ij})q}\\
&+&\frac{1}{p}F_{i}\sum\limits_{j=1}^{n}\xi_{j}|b_{ji}|^{\beta_{ji}p}e^{\epsilon
\tau_{ji}}+
\frac{1}{q}\xi_{i}\sum\limits_{j=1}^{n}F_{j}|b_{ij}|^{(1-\beta_{ij})q}
e^{\epsilon \tau_{ij}} \le 0 \label{Th61}
\end{eqnarray}
\begin{eqnarray}
&&\nonumber
(-d_{i}+\epsilon)\xi_{i}+G_{i}\bigg[\xi_{i}|a_{ii}|
+\frac{1}{p}\sum\limits_{j\ne i}\xi_{j}|a_{ji}|\bigg]
+\frac{1}{q}\xi_{i}\sum\limits_{j\ne
i}G_{j}|a_{ij}|\\
&+&\frac{1}{p}F_{i}\sum\limits_{j=1}^{n}\xi_{j}|b_{ji}|e^{\epsilon
\tau_{ji}}+
\frac{1}{q}\xi_{i}\sum\limits_{j=1}^{n}F_{j}|b_{ij}|
e^{\epsilon \tau_{ij}} \le 0 \label{Th62}
\end{eqnarray}
Then there exist constants $\theta_{i}$ such that
\begin{eqnarray}
&&
(-d_{i}+\epsilon)\theta_{i}+G_{i}\bigg[\theta_{i}|a_{ii}|
+\sum\limits_{j\ne i}\theta_{j}|a_{ji}|\bigg]
+F_{i}\sum\limits_{j=1}^{n}\theta_{j}|b_{ji}|e^{\epsilon
\tau_{ji}} \le 0 \label{Th63}
\end{eqnarray}

{\bf Theorem 7}\quad Suppose that $1\le p<\infty$. If there are
 positive constants $\epsilon$,  $\theta_{i}$, $i=1,2\cdots,n$, such that
\begin{eqnarray}
&&
(-d_{i}+\epsilon)\theta_{i}+G_{i}\bigg[\theta_{i}|a_{ii}|
+\sum\limits_{j\ne i}\theta_{j}|a_{ji}|\bigg]
+F_{i}\sum\limits_{j=1}^{n}\theta_{j}|b_{ji}|e^{\epsilon
\tau_{ji}} \le 0 \label{Th71}
\end{eqnarray}
Then we can find  real constants $\alpha^{*}_{ij}$,
$\beta^{*}_{ij}$, positive constants  $\xi_{i}$,
  $i=1,2\cdots,n$, such that
\begin{eqnarray}
&&\nonumber
(-d_{i}+\epsilon)\xi_{i}+G_{i}\bigg[\xi_{i}|a_{ii}|
+\frac{1}{p}\sum\limits_{j\ne
i}\xi_{j}|a_{ji}|^{\alpha^{*}_{ji}p}\bigg]
+\frac{1}{q}\xi_{i}\sum\limits_{j\ne
i}G_{j}|a_{ij}|^{(1-\alpha^{*}_{ij})q}\\
&+&\frac{1}{p}F_{i}\sum\limits_{j=1}^{n}\xi_{j}|b_{ji}|^{\beta^{*}_{ji}p}e^{\epsilon
\tau_{ji}}+ \frac{1}{q}\xi_{i}\sum\limits_{j=1}^{n}F_{j}|b_{ij}|^{(1-\beta^{*}_{ij})q}
e^{\epsilon \tau_{ij}} \le 0 \label{Th72}
\end{eqnarray}
hold. However, the following inequalities
\begin{eqnarray}
&&\nonumber
(-d_{i}+\epsilon)\xi_{i}+G_{i}\bigg[\xi_{i}|a_{ii}|
+\frac{1}{p}\sum\limits_{j\ne i}\xi_{j}|a_{ji}|\bigg]
+\frac{1}{q}\xi_{i}\sum\limits_{j\ne
i}G_{j}|a_{ij}|\\
&+&\frac{1}{p}F_{i}\sum\limits_{j=1}^{n}\xi_{j}|b_{ji}|e^{\epsilon
\tau_{ji}}+ \frac{1}{q}\xi_{i}\sum\limits_{j=1}^{n}F_{j}|b_{ij}|
e^{\epsilon \tau_{ij}} \le 0 \label{Th73}
\end{eqnarray}
are not true generally.

Following example verify previous theorems

\noindent{\bf Example 1} \quad
Let $d_{1}=2, ~ d_{2}=11, ~a_{11}=1, ~a_{12}=3,~ a_{21}=3,~ a_{22}=1,~
b_{ij}=0,~ G_{1}=G_{2}=F_{1}=F_{2}=1$.

\noindent i)\quad In case $p=1$.
The conditions ({\ref{Co23}})
in Corollary 2 are
\begin{eqnarray*}
-(d_{1}-|b_{11}|F_{1})\theta_{1}+\bigg[a_{11}G_{1}\theta_{1}
+\theta_{2}|a_{21}|G_{1}\bigg]+\theta_{2}|b_{21}|F_{1}\\
=-\theta_{1}+3\theta_{2}<0\\
-(d_{2}-|b_{22}|F_{2})\theta_{2}+\bigg[a_{22}G_{2}\theta_{2}
+\theta_{1}|a_{12}|G_{2}\bigg]+\theta_{1}|b_{12}|F_{2}\\
=3\theta_{1}-10\theta_{2}<0
\end{eqnarray*}
The solution to these two linear inequalities is
$3<\frac{\theta_{1}}{\theta_{2}}<\frac{10}{3}$.

\noindent ii)\quad
In case $p=2$. Pick any $\theta_{1},$ $\theta_{2}$ satisfying
$3<\frac{\theta_{1}}{\theta_{2}}<\frac{10}{3}$, for example,
 $\theta_{1}=19,\ \theta_{2}=6$
satisfy $3<\frac{\theta_{1}}{\theta_{2}}<\frac{10}{3}$.
Let $\eta_{1}=\theta_{1}^{2},\ \eta_{2}=\theta_{2}^{2}$ and
\begin{eqnarray}
\alpha_{12}=\frac{1}{2}\ln_{a_{12}}\bigg(a_{12}\eta_{1}^{\frac{1}{2}}\eta_{2}^{-\frac{1}{2}}\bigg)
            =\frac{1}{2}\ln_{3}\frac{19}{2}\\
\alpha_{21}=\frac{1}{2}\ln_{a_{21}}\bigg(a_{21}\eta_{2}^{\frac{1}{2}}\eta_{1}^{-\frac{1}{2}}\bigg)
            =\frac{1}{2}\ln_{3}\frac{18}{19}
\end{eqnarray}

The inequalities in (\ref{Co21}) become
\begin{eqnarray*}
-(d_{1}-|b_{11}|F_{1})\xi_{1}+\bigg[a_{11}G_{1}\xi_{1}
+\frac{1}{p}|a_{21}|^{\alpha_{21}p}G_{1}\xi_{2}\bigg]
    +\frac{1}{q}|a_{12}|^{(1-\alpha_{12})q}G_{2}\xi_{1}\\
    +\frac{1}{p}|b_{21}|^{\beta_{21}p}F_{1}\xi_{2}
    +\frac{1}{q}|b_{12}|^{(1-\beta_{12})q}F_{2}\xi_{1}
    =\frac{-10}{19}\xi_{1}+\frac{9}{19}\xi_{2}<0\\
-(d_{2}-|b_{22}|F_{2})\xi_{2}+\bigg[a_{22}G_{2}\xi_{2}
+\frac{1}{p}|a_{12}|^{\alpha_{12}p}G_{2}\xi_{1}\bigg]
    +\frac{1}{q}|a_{21}|^{(1-\alpha_{21})q}G_{1}\xi_{2}\\
    +\frac{1}{p}|b_{12}|^{\beta_{12}p}F_{2}\xi_{1}
    +\frac{1}{q}|b_{21}|^{(1-\beta_{21})q}F_{1}\xi_{2}
    =\frac{19}{4}\xi_{1}-\frac{-21}{4}\xi_{2}<0
\end{eqnarray*}
The solution is $\frac{9}{10}<\frac{\xi_{1}}{\xi_{2}}<\frac{21}{19}$.

Although the parameters $\alpha_{ij}$, $\beta_{ij}$ and/or
$\xi_{i}$ in (\ref{Co21}) do exist, it is very difficult to search
for these parameters directly by solving the nonlinear
inequalities (\ref{Co21}). Therefore, the criteria is of no great
significance in practice. Instead, the parameters in (\ref{Co23})
can be found easily by solving the linear inequalities
(\ref{Co23}).

Moreover, if we take $\alpha_{12}=\alpha_{21}=\beta_{12}
=\beta_{12}=\frac{1}{2}$.
The inequalities in (\ref{Co22}) are
\begin{eqnarray*}
-(d_{1}-|b_{11}|F_{1})\xi_{1}+\bigg[a_{11}G_{1}\xi_{1}
+\frac{1}{2}|a_{21}|G_{1}\xi_{2}\bigg]
+\frac{1}{2}|a_{12}|G_{2}\xi_{1}\\
+\frac{1}{2}|b_{21}|F_{1}\xi_{2}
+\frac{1}{2}|b_{12}|F_{2}\xi_{1}=\frac{1}{2}\xi_{1}+\frac{3}{2}\xi_{2}<0\\
-(d_{2}-|b_{22}|F_{2})\xi_{2}+\bigg[a_{22}G_{2}\xi_{2}
+\frac{1}{2}|a_{12}|G_{2}\xi_{1}\bigg]
+\frac{1}{2}|a_{21}|G_{1}\xi_{2}\\
+\frac{1}{2}|b_{12}|F_{2}\xi_{1}
+\frac{1}{2}|b_{21}|F_{1}\xi_{2}=\frac{3}{2}\xi_{1}-\frac{17}{2}\xi_{2}<0
\end{eqnarray*}
In this case, there is no solution. It means that this
criterion is less effective.

\section{Conclusions}

In this paper, we address criteria of global exponential stability  for delayed periodic
dynamical systems in
terms of various $L^{p}$ ($1\le p<\infty$) norms. A general
approach to investigate global exponential stability in terms of
$L^{p}$ ($1\le p<\infty$) norms and sufficient conditions are
given. Comparisons of various stability criteria are given.
More importantly, it is
pointed out that sufficient conditions in terms of $L^{1}$ norm
are enough in practice. The criteria in terms of
$L^{p}$ ($1< p<\infty$) norms are of no great significance.

\end{document}